\documentclass[11pt]{amsart}
\usepackage[numbers]{natbib}
\usepackage{a4,a4wide}
\usepackage{leqno}
\usepackage{graphics}
\usepackage{amssymb}
\usepackage{amsmath}
\usepackage{amsthm}
\usepackage[latin1]{inputenc}
\usepackage{euscript}
\usepackage{palatino}

\usepackage{color}

\newtheorem{theorem}{Theorem}[section]
\newtheorem{lemma}[theorem]{Lemma}

\theoremstyle{definition}

\newtheorem{corollary}[theorem]{Corollary}

\theoremstyle{remark}
\newtheorem{remark}[theorem]{Remark}

\numberwithin{equation}{section}

 \numberwithin{equation}{section}

\newcommand{\beq}{\begin{equation}}
\newcommand{\eeq}[1]{\label{#1}\end{equation}}

\def\O{{\Omega}}
\def\o{{\omega}}

\def\eps{{\epsilon}}

\def\l{{\mathcal{L}}}
\def\m{{\mathcal{M}}}

\def\r{{\mathcal{R}}}

\def\L{{\mathcal{L}}}

\def\S{{\mathcal{S}}}
\def\M{{\mathcal{M}}}

\def\R{{\mathbb{R}}}

\def\N{{\mathbb{N}}}

\newcommand{\opl}[1]{\L[{#1}]}
\newcommand{\opm}[1]{\M[{#1}]}

\newcommand{\opr}[1]{\r[{#1}]}

\newcommand{\dem}[1]{\vskip 0.2\baselineskip \noindent {\bf{#1}}\vskip 0.2\baselineskip }
\newcommand{\fdem}{\vskip 0.2 pt \hfill $\square$ }

\newenvironment{theo*}[2]{~\smallskip\newline \noindent\textit{\textbf{#1}}\;\newline \textit{#2}}

\newtheorem{theo}{\textbf{Theorem}}[section]

\newtheorem{cla}[theo]{\textbf{Claim}}

\def\tilde{\widetilde}
\def\ackno#1{\par\medskip{\small\textbf{Acknowledgements.}~#1}}

\title{Harnack type inequality for positive solution of some integral equation}

\author{J{\'e}r{\^o}me Coville\\
~\\
\textit{\tiny Max Planck Institute for Mathematics in the Sciences\\
Inselstrasse 22\\
D-04103 Leipzig\\
Germany}}
\address{Max Planck Institute for Mathematics in the Sciences\\
Inselstrasse 22\\
D-04103 Leipzig\\
Germany}

\curraddr{ INRA, Equipe BIOSP\\
Centre de Recherche d'Avignon\\ 
Domaine Saint Paul, Site Agroparc\\
84914 Avignon cedex 9, France
}

\email{jerome.coville@avignon.inra.fr}
\thanks{}


\begin{document}
\begin{abstract}
In this paper, we establish some Harnack type inequalities satisfied by positive solutions of nonlocal inhomogeneous equations arising in the description of various phenomena ranging from population dynamics to micro-magnetism.
For regular domains, we also derive an inequality up to the boundary. The main difficulty in such context lies in a precise control of the solutions outside a compact set and the existence of local uniform estimates. We overcome this problem by proving a contraction result which makes the $L^1$ norms of the solutions on two compact sets $\o_1\subset\subset\o_2$ equivalent. We also construct the principal positive eigenfunctions associated to particular nonlocal operators by using the corresponding Harnack type inequalities. 
\end{abstract}

\maketitle

\section{Introduction and Main results}
In this paper, we investigate the existence of a Harnack type inequality for positive continuous  solutions of
\begin{equation}
\opl{u}=0\label{h-eq-harmo}
\end{equation}
where the operator is defined by 
\begin{equation}
\opl{u}:=\int_{\O}k(x,y)u(y)\,dy -b(x)u, \label{h-eq-gene}
\end{equation}
 with $\O\subset \R^d$, $k \ge 0$  and $b(x)\in C(\O)$. Precise assumptions on $\O, k $ and $b$ will be given later on. 

Such  type of linear operator has been widely used to describe the dispersal of a population in its environment in the following sense. As stated in \cite{F1,F2,HMMV}, if $u(y,t)$ is
thought of as a density at a location $y$ at a time $t$ and $k(x,y)$
as the probability distribution of jumping from a location $y$ to
a location $x$, then the rate at which the individuals from all other
positions are arriving to the location $x$ is
$$\int_{\O} k(x,y)u(y,t)\,dy.$$
On the other hand, the rate at which the individuals are leaving the location $x$ to travel to all other places is $-b(x)u(x,t)$.
The operator $\l$ is called nonlocal since the behaviour of a  function at a given point depends on its values at points some distance away rather than just nearby ones.

In the past few years, there has been an intense interest to use such nonlocal operators to model problems in mathematical physics or in ecology, see, among other references, \cite{DGP,F2,HMMV,KM,Sch}. In particular, much attention has been devoted to the study of the nonlocal reaction diffusion equation 
\begin{equation}
\frac{\partial u}{\partial t}= \opl{u} +f(x,u(x,t)) \quad \text{ in }\quad \O\times\R^+, \label{h-eq-erdnl}
\end{equation}
where the usual elliptic diffusion operator is replaced by the  operator $\l$, see for example \cite{Bates1997,CCEM,HMMV,KM}. Problems related to \eqref{h-eq-erdnl} with a homogeneous nonlinearity have been widely treated in the literature when $b(x)\equiv 1$ and the kernel take the form of a convolution operator, i.e. $k(x,y)=J(x-y)$ with $J$ a probability density. For instance, the works in the references \cite{Bates1997,CovilleJuly2006,CDM2,CD2,DOPT1} are devoted to the study of  travelling front solutions existing for the problem \eqref{h-eq-erdnl} with a homogeneous bistable or monostable nonlinearity, while \cite{Bates1999a,CR} deal with the study of the steady state solutions of the problem \eqref{h-eq-erdnl} with a logistic type, bistable or power-like nonlinearity. The particular instance of the parabolic problem in $\R^n$ when $f = 0$ is considered in \cite{CCR,CERW}. In the heterogeneous case ( $b$ and $k$ general), few results  are known. This is  due mostly to the lack of compactness properties of $\l+\lambda$ or of its inverse. We quote \cite{CDM1,HMMV} which deal with a convolution kernel in a periodic environment with a monostable type nonlinearity and \cite{CCEM} in the case of a linear heterogeneous dispersal process in $\R$. 
\medskip

When the kernel $k(x,y)$ is smooth and compactly supported, it is easy to see that, when applied to a smooth function, the operator $\l$ can be rewritten into the following form 
$$\opl{u}:=\opm{u}+ \opr{u}$$ with $\m$ an elliptic operator 
\begin{equation}
\opm{u}:=a_{ij}(x)\partial_{ij}u +b_i(x)\partial_iu +c(x)u,\label{h-eq-ell}
\end{equation} and $\r$ an operator involving derivatives of higher order than in $\m$.
Indeed, we have 
$$
\opl{u}=\int_{\O}k(x,y)[u(y)-u(x)]\, dy -c(x)u,
$$
with $c(x):=b(x)-\int_{\O}k(x,y)dy$. Hence, setting $z=x-y$ and performing a formal Taylor expansion of $u$ in the integral, we obtain 
$$\int_{x-\O}k(x,x-z)[u(x-z)-u(x)]\, dy=a_{ij}(x)\partial_{ij}u +b_i(x) \partial_iu + R[x,\partial_{ijk}u]$$
with $a_{ij}(x)$ and $b_i(x)$ defined by the following expressions 
\begin{align*}
&a_{ij}(x)=\frac{1}{2}\int_{x-\O}k(x,x-z)z_iz_j\,dz\\
&b_{i}(x)=\int_{x-\O}k(x,x-z)z_i\,dz
\end{align*}
and $$ R[x,\partial_{ijk}u]:=\int_{x-\O} k(x,x-z)z_iz_jz_k\left(\int_0^{1}\int_0^{1}\int_0^{1}st^2 \partial_{ijk}u(x-hstz)\,dhdsdt\right)dz.$$
Therefore, for suitable kernels $k$, it seems reasonable to expect that the operator $\l$ shares many of the properties of $\m$. 
 
For a uniform elliptic operator $\m$, it is well known that the positive  solutions  of the equation 
$\opm {u}=0$ are satisfying an Harnack inequality, see \cite{E,GT}. That is to say
\begin{theorem}[Harnack Inequality \cite{E,GT}]\label{h-thm-harnackell}
Let $\O\subset \R^n$ be a domain and let $\o \subset \subset \O$ be a compact sub-domain. Then there exists a constant $C(\o)$ such that  for all positive smooth solutions  $u$ so that $$ \opm{u}=0$$  we have 
$$ u(x)\le C u(y) \quad \text{ for  all }\quad x,y \in \o.$$ 
\end{theorem}
Such estimate is an extremely important tool in the study of partial differential equations since they are providing in particular various key estimates in the analysis of the regularity of the solutions of PDE's. Moreover, such estimate plays a central role in the construction of a positive eigenfunction of $\m$, which is also an essential tool in the analysis of various nonlinear partial differential equations, see for example \cite{Berestycki2005a,BNV,BR}.

As observed in \cite{CCEM,CDM2,HMMV} for some particular choice of  $b$ and $k$, the principal eigenvalue $\lambda_1$ of $\l$ and its corresponding positive eigenfunction $\phi_1$, that is
\begin{equation}
\opl{\phi_1}+\lambda_1 \phi_1=0 \quad\text{ in }\quad \O, \label{h-eq-b}
\end{equation}
 are central in the analysis of the nonlocal reaction diffusion problems such as \eqref{h-eq-erdnl}. However,  due to the lack of compactness of $\l+\lambda$ or of its inverse,
 proving the existence of a principal eigenpair in this context is a very difficult task and no general result is known. As for elliptic operators, an Harnack type inequality for the positive  solutions of such nonlocal equation \eqref{h-eq-b} is expected to provide \textit{a priori} estimates for the construction of eigenpairs of a general operator $\l$. It is therefore of great interest to investigate the existence of such estimates. 
\medskip
 For some particular kernel $k(x,y)$, the Harnack estimate for solutions of \eqref{h-eq-harmo} are often obtained   as a consequence of the harmonic property of the solution. Indeed, when  $k(x,y)$ takes the following form  
$$k(x,y)=\begin{cases}
\frac{1}{|B(x,r(x))|} &\quad \text{ if } y  \in B(x,r(x))\\
 0 &\quad \text{ otherwise}
\end{cases}  
$$
where $c_0<r(x)<d(x,\complement\O)$ is a given function, then the positive solutions of \eqref{h-eq-harmo} are harmonic functions. There exists a vast literature  dealing with what is called in potential theory \textit{"The Converse Mean Value Problem"}. The central question is then to find the conditions on  $k(x,y)$ so that the solutions $u$ of \eqref{h-eq-harmo} are harmonic.  For more details on this subject, we refer the interested reader to \cite{Baxter1972,Cornea1995,Hansen1996,Hansen2008,Hansen2001,Hansen1993a,Netuka1994,VeechWinter1975}  and references therein. 
For the above kernel, to our knowledge, Veech \cite{VeechWinter1975} was the first to obtain a Harnack type estimate for the positive solutions of \eqref{h-eq-harmo} as a consequence of the restricted mean value property itself rather than the harmonic property of the solution. Later, Cornea and Vesel\'y \cite{Cornea1994,Cornea1995} have extended   the Harnack type estimate obtained by Veech to some more general kernels. More precisely, they have  obtained the following  

\begin{theorem}[ \cite{Cornea1994,Cornea1995}]
 Let  $k(x,y)$ be a kernel so that  $ \forall x\in \O$, there exists $  V_x$ and $W_x,$ two compact neighbourhoods of  $x$  and two strictly positive constants $ m_x$ and $M_x$     such that
  $W_x\subset \O$ and for any $y\in V_x$ we have
 \begin{equation}
 m_x \chi_{_{V_x}}\le \frac{1}{b(.)}k(.,y)\le M_{x} \chi_{_{W_x}}\label{h-cornea-cond}
 \end{equation}
where $\chi_{_{A}}$ is the characteristic function of the set $A$.
 Then, for all compact set $\o\subset \O$ there exists a constant $C(\o,k)$ such that for all super-median function $u$ (i.e. any function $u$ satisfying  $u(x)\ge \frac{1}{b(x)}\int_{\O}k(x,z)u(z)\, dz$ for all $x$) and any two points $x,y\in \o$, one have
 $$ \frac{1}{b(x)}\int_{\O}k(x,z)u(z)\,dz\le C\frac{1}{b(y)}\int_{\O}k(y,z)u(z)\,dz.$$ 
 \end{theorem}
Some Harnack inequality for nonlocal operators have also been proved for singular kernel, essentially for the generators of  pure jump processes, see for example \cite{Bass2004,Bass2005,Bass2002,Bog,Caffarelli2009,Kas}. In this context, the operator $\l$ takes the form 
\begin{equation}\label{h-eq-jp}
\opl{u}:=\int_{\O}k(x,y)[u(y)-u(x)]\,dy,
\end{equation}
where $k(x,y)\sim \frac{1}{|x-y|^{d+\alpha}}$ with $\alpha >0$. This particular structure implies some extra regularity for the solutions of \eqref{h-eq-harmo} and plays an essential role in the derivation of the Harnack inequality. 
 
In this work, we investigate the existence of a Harnack type inequality for a class of nonlocal operators $\l$ for which the  condition on the kernel \eqref{h-cornea-cond} is not always verified and when the compact set $\o$ can touch the boundary $\partial \O$ when it exists.  
More precisely, we study the existence of a Harnack type inequality for the positive  solutions of the operator defined in \eqref{h-eq-gene}, for any positive function $b(x)$ and when the kernel $k(x,y)$ takes the form introduced by Cortazar \textit{et al.} in \cite{CCEM}:
\begin{equation}
k(x,y) = J\left(\frac{x-y}{g(y)}\right)\frac{1}{g^n(y)},\label{h-eq-ccem}
\end{equation}
where $J$ is a probability density and the function $g$ is bounded and non negative.
We are particularly interested in finding some simple conditions on $J$, $g$, $b$ and $\O$ such that a Harnack type inequality holds for the positive solutions of the following equation 
\begin{equation}
\int_{\O}J\left(\frac{x-y}{g(y)}\right)\frac{u(y)}{g^d(y)}\, dy -b(x)u=0. \quad\text{ in }\quad \O. \label{h-eq}
\end{equation}

Now let us  state the precise assumptions that we impose on $J$, $g$ and $b$. Throughout this paper and without further notice we will always assume that $J$ $b$ and $g$ satisfy the following assumptions
\begin{align*}
&J \in L^{\infty}(\R^d)\cap L^{1}(\R^d), &\qquad(H1)\\
 & m_0\chi_{_{B(0,r_0)}}\le J\le M_0 \chi_{_{B(0,R_0)}}\; \text{for some positive constants }\, \ r_0, R_0, m_0,M_0, &\qquad(H2)\\
&g\in L^{\infty}(\O), 0\le g \le \beta,\, \frac{1}{g^d} \in L^p_{loc}(\O) \;\text{ with }\; p>1, &\qquad(H3)\\
&b \in C(\bar\O), b(x)>0. &\qquad(H4)
\end{align*}
Let us also denote by $\S$ the set of point where $g$ is vanishing (i.e. $\S:=\{x\in \bar\O|g(x)=0\}$ ). 
 Without loss of generality we will also assume that $R_0=1$.  Note that when  $\S\not =\emptyset$ the kernel $k(x,y)=J\left(\frac{x-y}{g(y)}\right)\frac{1}{g^d(y)} $ does not  satisfy automatically the assumption of Cornea and Vesel\'y (i.e. the condition \eqref{h-cornea-cond}) moreover the operator $\l$ can not be rewritten as a jump operator like in the equation \eqref{h-eq-jp}.

Under the above assumptions on $J$, $g$, and $b$, we first establish some kind of uniform estimate satisfied by the positive solutions of \eqref{h-eq}. More precisely,  we prove the following

\begin{theorem}\label{h-th2}
Let $J$, $g$, and $b$ satisfying (H1-H4) and  assume that $\O\cap\S\subset\subset \O$ and let $\o\subset\bar \O$ be a compact set. Let us  denote $\O(\o)$  the following set 
$$\O(\o):=\bigcup_{x\in \o}B(x,\beta).$$ 
Then there exists a positive constant $\eta^*$ such that for any $0<\eta\le\eta^*$ there exist a compact set $\o'\subset \subset \O(\o)\cap \O$ and a constant $C(J,\o,\o',b,g,\eta)$ such that the following assertions are verified
\begin{itemize} 
\item[(i)] $\{x\in \O(\o)\cap W_\eta|d(x,\partial(\O(\o)\cap W_\eta))>\eta\}\subset \o'$, where $W_\eta:=\{x\in \O|g(x)>\eta\}$
\item[(ii)] for all positive continuous solutions $u$ of \eqref{h-eq}, the following inequality holds: 
$$ u(x)\le C u(y) \quad \text{ for all }\quad x \in \o ,y\in \o'\cap \o .$$
\end{itemize}
\end{theorem}

Under an additional assumption on the regularity of the compact set $\o$, we have a more precise estimate. Namely, we prove the following result.
\begin{theorem}\label{h-th3}
Let $J$, $g$, and $b$ satisfying (H1-H4) and assume that $\O\cap\S\subset\subset \O$. Let $\o\subset\bar \O$ be a compact set which satisfies an uniform inner cone condition. Then there exists a positive constant $\eta^*$ such that for any $0<\eta\le\eta^*$, there exists a constant $C(J,\o,b,g,\eta)$ such that for all positive continuous solutions $u$ of \eqref{h-eq} the following inequality holds
$$ u(x)\le C u(y) \quad \text{ for all }\quad x \in \o ,y\in \o\cap \{y\in\O|g(y)>2\eta\}.$$
\end{theorem}

Note that in the above Theorems (Theorems  \ref{h-th2} and \ref{h-th3}), no condition is required on the open set $\O$. Therefore, the respective  inequalities are holding as well when $\O=\R^n$. In this particular case the assumption  on $\S$ can be weakened. More specifically instead of assuming that $\O\cap\S\subset\subset \O$ we can require that for any subset $\tilde\S$ of $\S$, there exists a ball $B(x_0,R)$ such that $$B(x_0,R)\cap\tilde\S\subset\subset B(x_0,R).$$ 

\begin{corollary}\label{h-cor3}
Let $J$, $g$, and $b$ satisfying (H1-H4) and assume that $\S$ satisfies the above condition and let $\o\subset\R^n$ be a compact set. Then there exists a positive constant $\eta^*$ such that for any $0<\eta\le\eta^*$ there exists a constant $C(J,\o,b,g,\eta)$ such that for all positive continuous solutions $u$ of \eqref{h-eq} the following inequality holds 
$$ u(x)\le C u(y) \quad \text{ for all }\quad x \in \o ,y\in \o\cap \{y\in\O|g(y)>2\eta\}.$$ 
\end{corollary}

Observe that  when the compact set $\o\subset \O$ and the considered kernel $k(x,y)$ satisfies the condition \eqref{h-cornea-cond} we recover the Harnack type estimate obtained by Cornea in \cite{Cornea1994,Cornea1995}.

As consequence of Theorems \ref{h-th2} and \ref{h-th3} we get various uniform estimates on the positive solutions $u$ of \eqref{h-eq} when $\O$ is a bounded domain. More precisely for general bounded domains we have  
 \begin{corollary}
Assume $J,g$ and $b$ satisfy (H1-H4),  let $\O$ be a bounded domain and assume that $\frac{1}{g^d}\in L^p(\O)$ for $p>1$. Then there exists a positive constant $\eta^*$ such that for any $0<\eta\le\eta^*$ there exists a compact set $\o'\subset\subset \O$ and a constant $C(J,\O,\o',b,g,\eta)$ such that the two following assertions hold 
\begin{itemize} 
\item[(i)]$\{x\in \O|d(x,\partial W_\eta)>\eta\}\subset \o'$
\item[(ii)] for all positive continuous  solutions $u$ of \eqref{h-eq}, 
$$ \sup_{\O}u\le C u(y) \quad \text{ for all } y\in \o'.$$
\end{itemize}
\end{corollary}
Note that in this case, the assumption $\S\cap \bar \O \subset\subset \O$ is not  any more required and the function $g$ can vanish at the boundary.    

In the case of a regular bounded domain and assuming that $g>\alpha$ the above inequality extends up to the boundary.
\begin{corollary}\label{h-cor2}
Let $J$, $g$, and $b$ satisfying (H1-H4) and assume that  $g\ge \alpha >0$ and $\O$ is a bounded domain which satisfies an uniform inner cone condition. Then there exists a constant $C(\O,g,J,b)$ such that for all positive continuous  solutions $u$ of \eqref{h-eq} it holds $$ \sup_{\O}u \le C \inf_{\O} u.$$
\end{corollary} 

Note that, in this particular case, whereas the classical Harnack type inequality remains true for positive solutions of the uniformly elliptic equation \eqref{h-eq-ell} its extension up to the boundary (Corollary \ref{h-cor2}) does not. The validity of such an extension is a consequence of the nonlocal nature of the equation considered.
\medskip

Our last result is an application of these estimates to the construction of a positive solution of \eqref{h-eq} for a particular $b(x)$. Let us consider the equation
\begin{equation}\label{h-eq-examp}
\int_{\O}J\left(\frac{x-y}{g(y)}\right)\frac{u(y)}{g^n(y)}\, dy - a(x)u =0 \quad\text{ in }\quad \O, 
\end{equation}
where $a(x)$ is defined as follows
\begin{equation*}
a(x):=\left\{
\begin{array}{l}
\int_{\O}J\left(\frac{y-x}{g(x)}\right)\frac{dy}{g^n(x)} \quad \text{ if } \quad x\not \in \S\\
1 \quad \text{ otherwise } 
\end{array}\right.
\end{equation*}
In the literature this equation corresponds to the nonlocal analogue of the usual eigenvalue problem for an elliptic operator with homogeneous Neumann boundary condition. For some particular situations, the existence of a positive solution is also well known. Indeed   when $g\equiv cste$ then  it is easy to see that for any $J$ any  constant is a solution of the problem \eqref{h-eq-examp}. In this  particular situation the kernel is of Markov type which is not true when $g\not  \equiv cste$.

Another well known case is when $\O$ is a bounded domain and $J$ is assumed to be a smooth (i.e. at least continuous ) function. Then in such a case the equation \eqref{h-eq-examp} reduces to a  homogeneous second kind Fredholm integral equation,  well studied in the literature, see for example \cite{Polyanin2008}.  
To our knowledge besides the two above situations the existence of a non trivial positive  solution of \eqref{h-eq-examp}  under  some general assumptions on $J,g$ and $b$ has only been obtained in the case $ \O=\R$ see  \cite{CCEM}.

Let us now state our result.
\begin{theorem}\label{h-th4}
Let $J$, $g$, and $b$ satisfying (H1-H4) and assume that $\O\cap\S\subset\subset \O$. Then there exists a positive continuous solution $p$ of \eqref{h-eq-examp}.
\end{theorem}

\medskip

\section{Comments and straightforward generalisations}
Before going into the proofs of these results, let us make some comments and explain our strategy of proofs.

We first want  to emphasize that the continuity assumption $H4$ made on  $b$ is not necessary and  all the results can be proved assuming that  $b\in L^{\infty}(\O)$ with $inf_{\O} b>0.$ We want also to  point out that all the above results  extend obviously to the kernels $k(x,y)$ of the form
$$k(x,y)=J\left(\frac{x_1-y_1}{g_1(y)};\frac{x_2-y_2}{g_2(y)};\ldots; \frac{x_n-y_n}{g_n(y)}\right)\frac{1}{\prod_{i=1}^ng_i(y)}.$$ 
We also remark that the assumption  $J>c_0\chi_{B(0,r_0)}$ cannot  be removed easily and a generalisation of the Harnack type estimate for more general measure $J$ seems delicate. Indeed  it is well known that the Harnack type inequality is false for some discrete Laplacian $\Delta_h$ which  corresponds to have in our framework  $b=g\equiv 1$ and $J:=\frac{1}{2}(\delta_{h}+\delta_{-h})$ where $\delta_h$ denotes the Dirac mass at the point $h$, see \cite{NetukaPrague1974}. For  example  by taking $h=2\pi$ we have $1+cos(x)$ is a non negative solution of the average equation $$u(x)=\int_{\R}J(x-y)u(y)\,dy=\frac{1}{2}(u(x+2\pi)+u(x-2\pi)). $$    
Extension of   these Harnack type inequality for positive measure $J$ are currently under consideration.

We want also to emphasize even though the average equation  \eqref{h-eq} have  some similarities with the average equations    
studied in the "Converse Restricted  Mean Value Problem", the two equations are fundamentally different. Indeed whereas a function satisfying a Restricted  Mean Value Property with respect to a measure $\mu$ verifies 
\begin{equation}\label{h-eq-mean}
 u(x)=\int_{B(x,\delta(x))}u(y)d\mu(y)\end{equation} for some  function $\delta(x)$,    a  solution $u$ of \eqref{h-eq} satisfies
$$ u(x)=\int_{B(x,g(y))\cap \O}u(y)d\nu(y),$$
 which in the case $g\not\equiv cste$ cannot be rewritten in the form \eqref{h-eq-mean} for some measure $\mu$.     

We also have recently observed that all our results on Harnack type inequality can be extended to the framework of super-median function (i.e $u$ so that $\opl{u}\le 0 $). In this framework, we have a Harnack type estimate  of the form
  \begin{theorem}\label{h-thgen1}
Let $J$, $g$, and $b$ satisfying (H1-H4).  Assume that  $\O\cap\S\subset\subset \O$ and   let $\o\subset\subset \O$ be a compact set. Let $\O(\o)$ denote the following set 
$$\O(\o):=\bigcup_{x\in \o}B(x,d_\o),$$ 
where $d_\o:=d(\o,\partial \O)$.
Then there exists a positive constant $\eta^*$ such that, for any $0<\eta\le\eta^*$, there exists a compact set $\o'\subset \subset \O(\o)$ and a constant $C(J,\o,\o',b,g,\eta)$ such that the following assertions are verified
\begin{itemize} 
\item[(i)] $\{x\in \O(\o)\cap W_\eta|d(x,\partial(\O(\o)\cap W_\eta))>\eta\}\subset \o'$, where $W_\eta:=\{x\in \O|g(x)>\eta\}$
\item[(ii)] for all positive continuous super median function $u$ of \eqref{h-eq} the following inequality holds: 
$$ \frac{1}{b(x)}\int_{\O}J\left(\frac{x-z}{g(z)}\right)\frac{u(z)\,dz}{g^d(z)}\le C \frac{1}{b(y)}\int_{\O}J\left(\frac{y-z}{g(z)}\right)\frac{u(z)\,dz}{g^d(z)}\quad \text{ for all }\quad x \in \o ,y\in \o'\cap \o .$$
\end{itemize}
\end{theorem}

We have  remarked that some of our  results can be generalized  easily to the situation where the function $g$ is vanishing at the boundary of the set $\O$ provided that $g(x)<d(x, \partial \O)$. More precisely, we have 
\begin{theorem}\label{h-thgen}
Let $J$, $g$, and $b$ satisfying (H1-H4).  Assume that  $\S=\partial \O \cup \S_i$ with $\O\cap\S_i\subset\subset \O$ and  $g(x)<d(x,\partial \O)$. Let $\o\subset\subset \O$ be a compact set and let $\O(\o)$ denote the following set 
$$\O(\o):=\bigcup_{x\in \o}B(x,d_\o),$$ 
where $d_\o:=d(\o,\partial \O)$.
Then there exists a positive constant $\eta^*$ such that, for any $0<\eta\le\eta^*$, there exist a compact set $\o'\subset \subset \O(\o)$ and a constant $C(J,\o,\o',b,g,\eta)$ such that the following assertions are verified
\begin{itemize} 
\item[(i)] $\{x\in \O(\o)\cap W_\eta|d(x,\partial(\O(\o)\cap W_\eta))>\eta\}\subset \o'$, where $W_\eta:=\{x\in \O|g(x)>\eta\}$
\item[(ii)] for all positive continuous solutions $u$ of \eqref{h-eq}, the following inequality holds: 
$$ u(x)\le C u(y) \quad \text{ for all }\quad x \in \o ,y\in \o'\cap \o .$$
\end{itemize}
\end{theorem}

Now let us  explain our  strategy to obtain such Harnack type inequality. It is mainly based on the following observation. For the harmonic functions (i.e functions $u$ such that $\Delta u=0$) it is well known (see \cite{E,GT})  that they satisfy a mean value equality 
$$ u(x)=\int_{B(x,r)}u(y)\,\frac{dy}{|B(x,r)|},$$ 
which holds for any ball $B(x,r)\subset\subset\O$. 
An Harnack inequality is derived  easily  from this property. Our idea is then to view a positive solution $u$ of \eqref{h-eq} as a positive function satisfying some mean value equality 
\begin{equation} u(x)=\frac{1}{b(x)}\int_{\O}u(y)\,d\mu(x,y), \label{h-eq-meq}\end{equation}
for some given measure $d\mu$, and to use this formulation to obtain some uniform estimates depending only on $\o$, $J$ and $b$. 
However, the later mean value property is fundamentally different from the one satisfied by harmonic functions in at least two ways. First, the measure $d\mu(x,y)$ is no more homogeneous and may be singular in the variable $y$. Second, the solution of equation \eqref{h-eq} satisfies the mean value equality for the fixed domain $\O$, whereas for harmonic functions the mean value equality holds for any ball compactly included in $\O$. All the difficulty in obtaining such Harnack type estimates arises from these two differences. 
We note that our proofs  differ from the one given in \cite{Cornea1994,Cornea1995} providing  alternative proofs of these Harnack type estimates.  In particular, from our analysis, we can extract an estimate of  the constant $C$ involved in the Harnack type inequality.

\medskip

 The paper is organized as follows. In the next section, we establish some general estimates that we use along our paper. Then we establish some uniform estimates satisfied by the positive solutions of \eqref{h-eq}. Next, in section \ref{h-s-h}, we prove the various Harnack type inequalities appearing in Theorems \ref{h-th2}, and \ref{h-th3}, and Corollaries \ref{h-cor3} to \ref{h-cor2}. Finally, the last section is devoted to the construction of a positive eigenfunction (Theorem \ref{h-th4}).

\section{Preliminaries}
In this section, we prove  some  technical lemmas  concerning some  sets of  positive functions satisfying a pointwise estimate. Similar estimates and inequalities can be found in  \cite{Cornea1995}.  Let us state our first lemma, 
\begin{lemma}\label{h-lem-tec}
Let $\O\subset \R^n$ be a domain and  let $x \in \O$ and $\eta>0$ so that  $B(x,4\eta)\subset \O$. Let $C_0$ be a positive constant and define $X$  the following set of functions 
 $$ X(x):=\left\{u\in C(\O)\;|\;  u\ge 0, \;  \text{ such that }\; u(y)\ge C_0 \int_{B(y,\eta)}u(s) \, ds \quad \text{for all }\; y\in B(x,2\eta)\right\}$$ 
Then for all  $\frac{\eta}{4}\le r\le \eta$, there exists a positive constant $C_1(r,\eta,C_0)\le 1$  so that for all  $u\in X$ we have
$$\int_{B(x,r)}u(s)\, ds\ge C_1 \int_{B(x,r+\frac{\eta}{4})}u(s)\,ds. $$  
\end{lemma}

Assume for the moment that Lemma \ref{h-lem-tec} holds. Then we can derive the following uniform estimates

\begin{lemma}\label{h-lem-loc-esti3}
Let $\O\subset \R^n$,  $\Sigma\subset\subset \O$ and $C_0>0$ be respectively a domain,  a smooth connected compact set and a positive constant.  Choose $\eta>0$  such that $$\O_{2\eta}:=\bigcup_{x\in \Sigma}B(x,4\eta) \subset\subset \O$$ and consider  $\tilde X$ the following set   of functions
$$ \tilde X:= \left\{ u \in C(\O)\,|\, u\ge 0, \; \text{ such that }\; u(y)\ge C_0 \int_{B(y,\eta)}u(s) \, ds \quad \text{for all }\; y\in \O_\eta:=\bigcup_{x\in \Sigma}B(x,2\eta)\right\} $$ 
Then there exists two positive constants $C(\eta,C_0,\Sigma)$ and $N(\eta, \Sigma)$ so that  
\begin{itemize}
\item[\it (i)]
for all $x,y\in \Sigma$ and for all $u\in \tilde X$ we have 
$$\int_{B(x,\frac{\eta}{4})}u(s)\, ds\ge C \int_{B(y,\frac{\eta}{4})}u(s)\,ds. $$ 
\item[\it (ii)] for all $x\in \Sigma$ and for all $u\in \tilde X$ we have 
$$\int_{B(x,\frac{\eta}{4})}u(s)\, ds\ge \frac{C}{N} \int_{\Sigma}u(s)\,ds. $$ 
\end{itemize}
\end{lemma}

\dem{Proof of Lemma \ref{h-lem-loc-esti3}:}
 
Let us start with the proof of \textit{(i)}.  Observe that for any $z\in \Sigma$ we have  by assumption 
\begin{equation}
B(z,2\eta)\subset \O_\eta  \text{ and  }  \tilde X \subset X(z).
\label{h-eq-loc-esti-1} 
\end{equation}
Therefore applying Lemma \ref{h-lem-tec} respectively  with $r=\frac{\eta}{4}$ and $ \frac{\eta}{2}$  yields for all $u \in \tilde X$  and for  all $z\in \Sigma$
\begin{align*}
& \int_{B(z,\frac{\eta}{4})}u(s)ds\ge C_1\left(\frac{\eta}{4},\eta, C_0\right)\int_{B(z,\frac{\eta}{2})}u(s)\, ds, \\
&\int_{B(z,\frac{\eta}{2})}u(s)ds\ge C_1\left(\frac{\eta}{2},\eta, C_0\right)\int_{B(z,\frac{3\eta}{4})}u(s)\, ds.
\end{align*}
Thus  for all $u \in  \tilde X$ and  for all $z\in \Sigma$ we have

\begin{equation}
  \int_{B(z,\frac{\eta}{4})}u(s)ds\ge C_2\int_{B(z,\frac{3\eta}{4})}u(s)\, ds, \label{h-eq-loc-esti2}
  \end{equation}
  with $C_2:=C_1\left(\frac{\eta}{4},\eta,C_0\right)C_1\left(\frac{\eta}{2},\eta,C_0\right)$.

  Using now that $\Sigma$ is a compact connected set then there exists a finite number of balls $ B(z_i,\frac{\eta}{4})$ covering $\Sigma$. That is to say, for some $N(\Sigma,\eta)\in\N$,
$$\Sigma\subset \bigcup_{i=1}^{N}B(z_i,\frac{\eta}{4}).$$

Now let $x,y \in \Sigma$ be fixed . From the covering,  we  can  find a finite sequence   $(t_n)_{n\in \{0,\ldots,N_0\}}$ of elements of  $\Sigma$   so that
$$
\begin{cases}
&N_0\le N(\Sigma,\eta)+2\\
&t_i \in \{x,y,z_1, \ldots ,z_N \},\, t_0=x \quad\text{ and }\quad t_{N_0}=y\\ 
& \forall i\in \{0,\ldots, N_0-1\}, B(t_{i+1},\frac{\eta}{4})\subset B(t_i,\frac{3\eta}{4})
\end{cases}
$$ 
Since for all $u\in \tilde X$, $u\ge 0$, using  \eqref{h-eq-loc-esti2} with $z=x$ and the definition of $t_1$, we deduce that for all $u \in \tilde X$  
$$ \int_{B(x,\frac{\eta}{4})}u(s)\,ds\ge C_2\int_{B(x,\frac{3\eta}{4})}u(s)\, ds \ge C_2\int_{B(t_1,\frac{\eta}{4})}u(s)\,ds. $$
Then by induction,   for all $u \in \tilde X$ we achieve
$$  \int_{B(x,\frac{\eta}{4})}u(s)\,ds\ge (C_2)^{N_0} \int_{B(y,\frac{\eta}{4})}u(s)\,ds.$$
Since  $N(\eta,\Sigma)$ is independent of $x$ and $y$  and $N_0\le N+2$ we have 
for all $x,y \in \Sigma$ and for all $u\in \tilde X$,
  
$$  \int_{B(x,\frac{\eta}{4})}u(s)\,ds\ge (C_2)^{N+2} \int_{B(y,\frac{\eta}{4})}u(s)\,ds,$$
which proves \textit{(i)}.

To obtain \textit{(ii)}, we just have to remark that using the above covering of $\Sigma$  and  \textit{(i)} for all $x\in \Sigma$ and all $u \in \tilde X$ we have
\begin{align*}
 \int_{B(x,\frac{\eta}{4})}u(s)\,ds&\ge \sum_{i=1}^{N} \frac{(C_2)^{N+2}}{N} \int_{B(z_i,\eta)}u(s)\,ds\\
&\ge \frac{(C_2)^{N+2}}{N} \int_{\bigcup_{i=1}^NB(y_i,\eta)}u(s)\,ds\\
&\ge  \frac{(C_2)^{N+2}}{N}\int_{\Sigma}u(s)\,ds.\\
\end{align*}
\fdem

Let us now turn our attention to the proof of the technical Lemma \ref{h-lem-tec}.
\dem{Proof of Lemma \ref{h-lem-tec}}

Since $r\le\eta$, for all $z\in B(x,r)$ we have  $B(z,\eta)\subset B(x,2\eta)$ and by assumption for all $u \in X(x)$ we have
 \begin{equation}
\int_{B(x,r)}u(s)\, ds \ge C_0\int_{B(x,r)}\left(\int_{B(z,\eta)}u(s)\, ds\right)\, dz.
\label{h-eq-tec2}
\end{equation}
Let us now consider the annulus $A:=A(x,r',r)$, for some $0<r'< r$ which will be chosen later on.
 Observe that for $\eps>0$, $A$ can be covered by a finite numbers of balls $B(z,r-r'+\eps)$, where $z\in \partial B(x,r')$.
Namely, we have for some $N_{(r-r'+\eps)}\in\N$, $$A\subset \bigcup_{i=1}^{N_{(r-r'+\eps)}}B(z_i,(r-r'+\eps)). $$
 Observe that by construction, for any    $N_{(r-r')+\eps)}-$tuplet $(\tilde z_1, \tilde z_2, \ldots, \tilde z_{N_{(n-r')+\eps)}})$ such that $ \tilde z_i\in B(z_i,(r-r')+\eps)\cap A$ we have 
 $$ A\subset \bigcup_{i=1}^{N_{(r-r'+\eps)}}B(\tilde z_i,2(r-r')+2\eps).$$

 Moreover, for $r-r'+\eps<\frac{\eta}{2}$ and  for any    $N_{(r-r')+\eps)}-$tuplet $(\tilde z_1, \tilde z_2, \ldots, \tilde z_{N_{(r-r')+\eps)}})$ such that $ \tilde z_i\in B(z_i,(r-r')+\eps)\cap A$, we see that 
  
 \begin{equation}
A(x,r',r'+\frac{\eta}{2}-[(r-r')+\eps])\subset \bigcup_{i=1}^{N_{(r-r'+\eps)}}B(\tilde z_i,\eta).\label{h-eq-geocond}
\end{equation}
 Indeed, let    $(\tilde z_1, \tilde z_2, \ldots, \tilde z_{N_{(n-r')+\eps)}})$ be a $N_{(r-r')+\eps)}-$tuplet such that $ \tilde z_i\in B(z_i,(r-r')+\eps)\cap A$ and take $y \in A(x,r',r'+\frac{\eta}{2}-[(r-r')+\eps])$. Then $d(y,\partial B(x,r'))\le \frac{\eta}{2}-[(r-r')+\eps]$ and there exists  $\bar z \in \partial B(x,r') $ so that $$\|y-\bar z\|=d(y,\partial B(x,r')).$$ 
 
 Since $\bar z \in \partial B(x,r')$, by construction there exists $z_j \in \partial B(x,r')$ so that $$\bar z \in B(z_j,r-r'+\eps).$$
 Choose  $\tilde z_j$ the corresponding element in the $N_{(r-r')+\eps)}-$tuplet $(\tilde z_1, \tilde z_2, \ldots, \tilde z_{N_{(n-r')+\eps)}})$ and compute $\|y-\tilde z_j\|$ then we have
 \begin{align*}
 \|y-\tilde z_j\|&\le\|y-\bar z\|+\|\bar z-z_j\|+\| z_j-\tilde z_j\|  \\
 &\le \frac{\eta}{2}+[(r-r')+\eps]\\
 &\le \eta
 \end{align*}
  since $[(r-r')+\eps]<\frac{\eta}{2}$. Thus $y\in B(\tilde z_j,\eta)$ and \eqref{h-eq-geocond} holds.
 
Let us fix $\eps=\frac{r}{8}$ and  choose $r'=\frac{15r}{16}$. Let us also denote $\mu(S)$ the Lebesgue measure of a set $S$ and consider $A_i:=B(z_i,(r-r')+\eps)\cap A$ . By construction, since $y_i\in \partial B(x,r')$, we have $\mu(A_i)=\mu(A_j)$ for all $i,j$ and for each $A_i$ from \eqref{h-eq-tec2} we have for all $u \in X(x)$
 
\begin{equation}\label{h-eq-3}
\int_{B(x,r)}u(s)\, ds \ge C_0\int_{A_i}\left(\int_{B(z,\eta)}u(s)\, ds\right)\, dz.
\end{equation}
Therefore for all $u\in X(x)$, on each $A_i$ there exists a point $\tilde z_i\in A_i$ which can depend on $u$  such that  
\begin{equation}
\int_{B(x,r)}u(s)\, ds\ge \frac{C_0}{\mu(A_i)}\int_{B(\tilde z_i,\eta)}u(s)\, ds. 
\end{equation}
Using that $\mu(A_i)=\mu(A_j)$ for all $i,j$, we deduce that for each $u \in X(x)$ there exists a $N_{(r-r')+\eps)}-$tuplet $(\tilde z_1, \tilde z_2, \ldots, \tilde z_{N_{(n-r')+\eps)}})$ so that  
 \begin{align*}
\int_{B(x,r)}u(s)\, ds&\ge \frac{C_0}{\mu(A_i)N_{(r-r'+\eps)}}\sum_{i=1}^{N_{(r-r'+\eps)}}\int_{B(\tilde z_i,\eta)}u(s)\, ds \\
&\ge\frac{C_0}{\mu(A_i)N_{(r-r'+\eps)}} \int_{\bigcup_{i=1}^{N_{(r-r'+\eps)}} B(\tilde z_i,\eta)}u(s)\, ds. 
\end{align*}

Since by construction for each $u\in X(x)$ we have  $\tilde z_i\in B(z_i,(r-r')+\eps)\cap A$ for all $i$.  Using the geometric condition \eqref{h-eq-geocond}, it follows that for all $u \in X(x)$ 
\begin{equation*}
\int_{B(x,r)}u(s)\, ds\ge \frac{C_0}{\mu(A_i)N_{r-r'+\eps}}\int_{A( x, r',r'+\frac{\eta}{2}-[(r-r')+\eps])}u(s)\, ds. 
\end{equation*} 
Therefore, 
\begin{align*}
\int_{B(x,r)}u(s)\, ds&\ge\frac{C_0}{\mu(A_i)2N_{r-r'+\eps}}\int_{A( x, r',r'+\frac{\eta}{2}-[(r-r')+\eps])}u(s)\, ds +\frac{1}{2}\int_{B(x,r)}u(s)\, ds \\
&\ge C_1\int_{B(x,r)\cup A( x, r',r'+\frac{\eta}{2}-[(r-r')+\eps])}u(s)\, ds\\
&\ge C_1\int_{B(x,r'+\frac{\eta}{2}-[(r-r')+\eps])}u(s)\, ds,\\
&\ge C_1\int_{B(x,r +h )}u(s)\, ds,
\end{align*}
where $C_1:=\min\{\frac{C_0}{2\mu(A_i)N_{r-r'+\eps}},\frac{1}{2}\}$ and  $h:=\frac{\eta}{2}- [2(r-r')+\eps]$.

Since  by construction $h=\frac{\eta}{2}- [2(r-r')+\eps]\ge \frac{\eta}{2}-\frac{r}{4}\ge \frac{\eta}{4}$ we achieve for all $u \in \tilde X(x)$,
$$
\int_{B(x,r)}u(s)\, ds\ge C_1\int_{B(x,r+h)}u(s)\, ds \ge C_1\int_{B(x,r+\frac{\eta}{4})}u(s)\, ds.
$$

\fdem

\begin{remark}
Note that from our construction we can make the constant $C_1$ independent of $r$. Namely, by our choice of $r'$ and $\eps$ $N_{r-r'+\eps}=N_1$ is invariant with $r$ and  we have $$C_1\ge  \min\{\frac{C_0}{2\eta^d\mu( A(0,\frac{15}{16},1)\cap B(z, \frac{17}{16}))N_1},\frac{1}{2}\},$$
where $z$ is any point of $\partial B(0, \frac{15}{16})$ and $d$ the dimension of space.
\end{remark}

\section{Local uniform estimates \label{h-s-lue}}
In this section, we establish some local uniform estimates, which will play an essential role in deriving Harnack type inequalities.
\begin{lemma}\label{h-lem-loc-esti1}
Let $\O\subset\R^n$ and $u\in C(\O)$ be respectively a connected domain and a positive solution of \eqref{h-eq}. Let $\O'\subset\O$  be a compact set such that $g\ge\alpha>0$ in $\O'$. Then 
there exists $\eps^*>0$ such that for all $\eps\le \eps^*$, there exists $\O_\eps$ and $C_\eps(\alpha,\beta,J,\eps,b)$ such that
$$\int_{\O_{\eps}}u(y)\,dy \ge C_\eps\int_{\O'}u(y)\,dy.$$
Moreover, $\O_\eps$ satisfies the following chain of inclusion 
$$\left\{x\in\O'| d(x,\partial\O')>\alpha\eps \right\}\subset\O_\eps\subset\left\{x\in\O'| d(x,\partial\O')>\frac{\alpha\eps}{2}\right\}.$$
\end{lemma}
\dem{Proof :}
Since $u$ is a positive solution of \eqref{h-eq}, using that $\O'\subset \O$, we have 
\begin{equation}\label{h-eq-super}
 \int_{\O'}J\left(\frac{x-y}{g(y)}\right)\frac{u(y)}{g^d(y)}\, dy-b(x)u(x)\le 0\qquad \text { in }\qquad \O.
\end{equation} 
The domain $\O'$ being compact and $u$ continuous, we can integrate \eqref{h-eq-super} over $\O_\eps\subset \subset \O'$ and we have
$$\int_{\O_\eps}\int_{\O'}J\left(\frac{x-y}{g(y)}\right)\frac{u(y)}{g^d(y)}\, dy \le \int_{\O_\eps}b(x)u(x)\,dx.$$
Using that $g\ge \alpha >0$ in $\Omega'$, Fubini's Theorem and setting $z=\frac{x-y}{g(y)}$ , we end up with
\begin{align}
 \int_{\O_\eps}b(x)u(x)\,dx&\ge\int_{\O'}\frac{u(y)}{g^d(y)}\left( \int_{\O_{\eps}}J\left(\frac{x-y}{g(y)}\right)\, dx\right)\, dy \\
&\ge\int_{\O'}u(y)\left( \int_{\O_{\eps,y}}J(z)\, dz\right)\, dy 
\end{align}
where $\O_{\eps,y}:=\frac{\O_\eps-y}{g(y)}$.
We claim that
\begin{cla}
There exist $\O_\eps$ and $c_0>0$ such that, for all $y\in \O'$, $$\int_{\O_{\eps,y}}J(z)dz>c_0.$$
\end{cla}
Observe that by proving the above claim the proof of the lemma is ended. Indeed, assuming it is true, then we derive from the above inequality
 \begin{align*}
 \int_{\O_\eps}b(x)u(x)\,dx&\ge\int_{\O'}u(y)\left( \int_{\O_{\eps,y}}J(z)\, dz\right)\,dy,\\
&\ge c_0\int_{\O'}u(y)dy.
\end{align*}
Hence, $$\int_{\O_\eps}u(x)\,dx\ge \frac{c_0}{\|b\|_{\infty}}\int_{\O'}u(y)dy.$$

\fdem

\dem{Proof of the claim}
By assumption, since $J(0)>0$ and $J$ smooth, there exist $r_0>0$ and $c_0$ such that $\min_{B(0,r_0)}J>c_0$.\\
Fix $\eps\le \min\{\frac{r_0}{2};\frac{1}{4}\}$ such that the two sets 
\begin{align}
\O'_\eps:=\{x\in \O'|\; d(x,\partial \O')\ge \eps\alpha \}\\
\tilde\O'_\eps:=\{x\in \bar\O'|\; d(x,\partial \O')\le \frac{\eps\alpha}{2} \}
\end{align}
are non empty disjoint sets.
Choose $\O_\eps$ smooth (at least $C^2$) so that $\O'_\eps\subset\subset\O_\eps$ and $\O_\eps\cap \tilde \O'_\eps=\emptyset$.
By construction, we see that  $\bar \O_\eps$ is compact and  for all $y\in \O'$,   $d(y,\O_\eps)<\eps\alpha$. 

Again  by construction we observe that  for  $\delta\le \frac{\eps\alpha}{4}$, we  also have    
$$\forall\; z\in \O_\eps,\quad B(z,\delta)\subset\O'.$$ 

From the uniform regularity of  $\O_\eps$, there exists  a constant $\delta$ small enough,  say  $0<\delta<\delta_1\le \frac{\eps\alpha}{4}$ where $\delta_1$ only depends on the regularity of $\partial \O_\eps$, such that for all $z\in \O_\eps$, there exists $ z' \in \O_\eps\cap B(z,\delta)$ satisfying 
\begin{equation}\label{h-eq-cla}
 B(z',\frac{\delta}{8})\subset B(z,\delta)\cap \O_\eps.
\end{equation}

Now, pick $y \in \O'$. Since $\bar\O_\eps$ is compact, there exists $z_0\in \bar\O_\eps$ such that $\|y-z_0\|=d(y,\O_\eps)$. 
Using \eqref{h-eq-cla}, it follows that there exists $z_0'$ such that 
\begin{equation}
  B(z_0',\frac{\delta}{8})\subset B(z_0,\delta)\cap \O_\eps. \label{h-eq-boule}
  \end{equation}
Recall that $\O_{\eps,y}= \frac{\O_\eps -y}{g(y)}$ thus from \eqref{h-eq-boule} it follows that 
$$\frac{B(z_0',\frac{\delta}{8})-y}{g(y)} \subset\O_{\eps,y}.$$
 
Take now $s\in B(z_0',\frac{\delta}{8}) $ and let us compute $ \frac{\|s-y\|}{g(y)}$: 
\begin{align*}
\frac{\|s-y\|}{g(y)}&\le \frac{\|s-z'_0\|+\|z_0-z'_0\|+\|z_0-y\|}{\alpha}\\
&\le\frac{\delta}{8\alpha} +\frac{\delta}{\alpha}+\eps.
\end{align*}
Since $\delta\le \frac{\eps\alpha}{4}$ and $\eps\le \frac{r_0}{2}$, 
 we achieve $$ \frac{\|s-y\|}{g(y)}\le ( \frac{1}{64}+\frac{1}{8} +\frac{1}{2})r_0 \le r_0.$$

Finally, let us compute $\int_{\O_{\eps,y}}J(z)dz$. From the above construction, we have 
\begin{align*}
\int_{\O_{\eps,y}}J(z)dz&\ge \int_{\frac{B(z_0',\frac{\delta}{8})-y}{g(y)}}J(z)\,dz\\
&\ge c_0 \int_{\frac{B(z_0',\frac{\delta}{8})-y}{g(y)}}dz\\
&\ge c_0 \mu(B(0,\frac{\delta}{8\beta})).
\end{align*}
Since the above computation is independent of $y\in\O'$, the claim is proved. 

\fdem
\begin{remark}
Observe that from the above computation, we have a certain degree of freedom over the parameter $\eps$, which  later can  be chosen at our convenience. 
\end{remark}

Let us now show another important estimate.

\begin{lemma}\label{h-lem-loc-esti2}
Let $\O\subset\R^n$ be a connected set and $u$ a positive continuous function satisfying \eqref{h-eq}. Let $\O'\subset\O$ be such that $g\ge \alpha>0$ in $\O'$. Then for any $\O''\subset\subset \O'$ there exist $\delta$ and a constant $C(\beta,J,\alpha,\delta,b)$ such that 
$$\forall\, x\in \O_\delta:=\bigcup_{x\in\O''}B(x,\delta),\quad u(x)\ge C\int_{B(x,\delta)}u(y)\,dy. $$ 
\end{lemma}

\dem{Proof:}
Let $d:=d(\O'',\partial \O')$. By assumption, one has $d>0$. 
Since $u$ is positive and $\O'\subset\O$, we deduce from \eqref{h-eq} that at $x\in\O$ we have
\begin{align*}
b(x)u(x)&\ge \int_{\O'}J\left(\frac{x-y}{g(y)}\right)\frac{u(y)}{g^d(y)}\, dy\\
&\ge\int_{B(x,\beta)\cap\O'}J\left(\frac{x-y}{g(y)}\right)\frac{u(y)}{g^d(y)}\, dy. \\
\end{align*}
Using that  $J(0)>0$, $\alpha\le g\le\beta$ in $\O'$ and $\| b\|_{\infty}<C$ we see that for $\delta$ small, say $\delta\le \delta_1$, we have 
 \begin{align*}
u(x)&\ge\int_{B(x,\delta)\cap\O'}J\left(\frac{x-y}{g(y)}\right)\frac{u(y)}{g^d(y)}\, dy\\
& \ge \frac{\min_{B(0,\delta)} J}{\beta^n\| b\|_{\infty}}\int_{B_\delta(x)\cap\O'}u(y)\, dy.
\end{align*}
Choosing $\delta<\min\{\frac{d}{2},\delta_1\}$, it follows that for any $x\in \O_\delta$, $B(x,\delta)\subset \O'$. Hence, for all $x\in \O_\delta$, we have 
$$ u(x)\ge \frac{\min_{B(0,\delta)} J}{\beta^n\| b\|_{\infty}}\int_{B_\delta(x)}u(y)\, dy$$
and the Lemma is proved.
\fdem

\section{Harnack type inequalities \label{h-s-h}}
We are now in position to prove the different Harnack type  inequalities, Theorems  \ref{h-th2} and \ref{h-th3} and Corollaries  \ref{h-cor2} and \ref{h-cor3}. 
The proof of the Corollaries  come as a straightforward application of the main Theorems and thus left to the reader. Simple proofs of Theorem \ref{h-th2} can be obtained using Theorems \ref{h-th3}, so let us first prove Theorem \ref{h-th3}.

\dem{Proof Theorem \ref{h-th3}:}
Before we begin, let us make some remarks and introduce some notation. 
First observe that if the estimates in Theorem \ref{h-th3} hold for a given compact set $\o\subset \O$, 
then they hold as well for any compact set $\tilde\o \subset\o$. Indeed, since the estimates in Theorem \ref{h-th3} hold for $\o$, there exists a positive constant $\eta^*(\o)$ so that for any $0<\eta\le\eta^*(\o)$ there exists a constant $C(\eta)$  such that  for all positive solution $u$ of \eqref{h-eq} the following inequality holds: 
$$ u(x)\le C u(y) \quad \text{ for all }\quad x \in \o ,y\in \o\cap W_{2\eta},$$
where $W_\eta:=\{y\in\O | g(y)\ge \eta\}.$

Using now  that any positive solutions of \eqref{h-eq} satisfies
$$\sup_{\tilde\o}u\le \sup_{\o}u\le C(\eta) \inf_{W_{2\eta}\cap \o}u\le C(\eta) \inf_{W_{2\eta}\cap \tilde\o}u,$$
it follows that  the estimates in Theorem \ref{h-th3} hold for $\tilde\o$.

From the above observation, we can restrict our attention to compact set $\o\subset \O$ such that $\S\subset\subset \o$. 
Fix now $\o$ and let us define the following sets 
\begin{align*}
&\o_\eta:=\bigcup_{x\in \o}B(x,\eta)\cap \O\\
&Z_\eta:=\{y\in\O | g(y)< \eta\}\\
&W_\eta:=\{y\in\O | g(y)\ge \eta\}.
\end{align*}
Since $\frac{1}{g^d}\in L^p_{loc}$ with $p>1$, we can  choose $\eta_1$ small enough such that
\begin{equation}
\int_{\o\cap Z_{\eta^*}}\frac{dy}{g^d(y)}\le \frac{\inf_{\o}b}{2\|J\|_{\infty}}.\label{h-eq-sing-regeta}
\end{equation}
Since $\S\subset\subset \o$, we can choose $\eta_1$ smaller if necessary to achieve $\o_{\eta_1}\cap Z_\eta\subset \o$.
Fix now, $0<\eta\le \eta_1$.
We are now in a position to prove the Theorem. The proof follows essentially four steps. 
\subsection*{Step 1:}
Now, define the following bounded set $$\O(\o):=\bigcup_{x\in \o}B(x,\beta),$$
and set the measure $d\mu=\frac{dy}{g^d(y)}$, which is well defined since $\frac{1}{g^d}\in L^1_{loc}$.
Since $J$ is compactly supported, it follows that in $\o$, $u$ satisfies 
\begin{align}
u(x) &=\frac{1}{b(x)}\int_{\O(\o)\cap\O} J\left(\frac{x-y}{g(y)}\right) u(y) d\mu(y)\\
 &=\frac{1}{b(x)}\int_{\O(\o)\cap Z_\eta} J\left(\frac{x-y}{g(y)}\right) u(y) d\mu(y)+ \frac{1}{b(x)}\int_{\O(\o)\cap W_\eta} J\left(\frac{x-y}{g(y)}\right) u(y) d\mu(y).\label{h-eq-sing-reg1}
\end{align}

Observe that for $x\in \o$, $y\in Z_\eta \cap(\O(\o)\setminus\o_\eta)$, we have 
$$\left|\frac{x-y}{g(y)}\right|\ge 1.$$ 
Therefore since $supp(J)\subset B(0,1)$, it follows that for $x\in \o$ 
$$ \frac{1}{b(x)}\int_{\O(\o)\cap Z_\eta} J\left(\frac{x-y}{g(y)}\right) u(y) d\mu(y)= \frac{1}{b(x)}\int_{\o_\eta\cap Z_\eta} J\left(\frac{x-y}{g(y)}\right) u(y) d\mu(y)$$ 
and from \eqref{h-eq-sing-reg1} we get 
\begin{equation}
u(x)\le \frac{1}{b(x)}\int_{\o_\eta\cap Z_\eta} J\left(\frac{x-y}{g(y)}\right) u(y) d\mu(y)+\frac{\|J\|_{\infty}}{\inf_{\o}b(x)}\int_{\O(\o)\cap W_\eta }u(y) d\mu(y). \label{h-eq-sing-reg2}
\end{equation}
Since $u$ is continuous and $\o$ is compact, $u$ achieves its maximum at some point, say $x_0\in\o$.
At this point, from \eqref{h-eq-sing-reg2} we have:
\begin{equation}
u(x_0)\le \frac{1}{b(x_0)}\int_{\o_\eta\cap Z_\eta} J\left(\frac{x_0-y}{g(y)}\right) u(y) d\mu(y)+\frac{\|J\|_{\infty}}{\inf_{\o}b(x)}\int_{\O(\o)\cap W_\eta }u(y) d\mu(y).
\end{equation}
Using that $\o_\eta\cap Z_\eta\subset \o$ and \eqref{h-eq-sing-regeta}, it follows that 
\begin{equation}
u(x_0)\le \frac{u(x_0)}{2}+\frac{\|J\|_{\infty}}{\inf_{\o}b(x)}\int_{\O(\o)\cap W_\eta }u(y) d\mu(y).
\end{equation}
Therefore, 
\begin{equation}
u(x_0)\le \frac{2\|J\|_{\infty}}{\inf_{\o}b(x)}\int_{\O(\o)\cap W_\eta }u(y) d\mu(y) \label{h-eq-sing-reg3}.
\end{equation}

\subsection*{Step 2:}
For any $\nu\in \R$, let us consider the set $\o_\nu:=\{x\in \o|d(x,\partial \o)\ge \nu\}$ and $C(x,\theta,a)$ be the cone issued from $x$ with angle $\theta$ and height $a$. 
On one hand, since $\S\subset\subset \o$, there exists $\nu_0>0$ such that $\o\setminus \o_{4\nu_0}\cap \S =\emptyset$.
On the other hand, since $\o$ satisfies an uniform inner cone condition, it follows that for $\nu$ small enough, say $\nu \le \nu^*$, there exists $r(\nu)>0$ such that for any $x\in \o\setminus \o_{\nu}$, there exists $\bar x \in \o$ such that 
\begin{align*}
& B(\bar x,r)\subset C_{x,\theta,a}\cap (\o_{\nu}\setminus \o_{4\nu})\\
&B(\bar x,r)\subset B(x,\beta).
\end{align*}
Let us now fix $\nu\le \min\{\nu_0,\nu^*\}$ and take $\eta^*:=\min_{\O(\o)\setminus \o_{4\nu}}g$. By construction, $\eta^*>0$.

Now take any $x \in \o\setminus \o_\nu$. From \eqref{h-eq}, using the uniform inner cone property, we have 
\begin{align}
u(x)&=\frac{1}{b(x)}\int_{\O}J\left(\frac{x-y}{g(y)}\right)\frac{u(y)}{g^d(y)}\, dy\\
&\ge \frac{1}{b(x)}\int_{C_{x,\theta,a}\cap B(x,\beta)}J\left(\frac{x-y}{g(y)}\right)\frac{u(y)}{g^d(y)}\, dy\\
&\ge \frac{1}{b(x)}\int_{B(\bar x,r)}J\left(\frac{x-y}{g(y)}\right)\frac{u(y)}{g^d(y)}\, dy.
\end{align}

Recall that $g\ge \eta^*$ in $B(\bar x, r)$. Therefore, since $J(0)>0$, $b>0$, there exists $\delta_0$ and $C_0$ independent of $x$ such that $B(\bar x,\delta_0)\subset B(\bar x,r)$ and 
 \begin{equation}
u(x)\ge C_0\int_{B(\bar x,\delta_0)}u(y)\, dy. \label{h-eq-sing-reg4}
\end{equation}

Fix now $\eta \le \min\{\frac{\eta_1}{2},\frac{\eta^*}{2}\}$ such that 
$W_\eta\setminus W_2\eta \subset\subset \o_\nu$, and let 
$$d:=d\Big(\o_\nu\cap W_{2\eta}, \partial(\O(\o)\cap W_{\eta})\Big).$$
By construction, we have $d>0$. Indeed, since $\eta \le \eta^*$, we have $\partial(\O(\o)\cap W_{\eta})=\Gamma_1 \cup \Gamma_2$ where $\Gamma_1 \subset \overline{(\O(\o)\setminus \o_{\frac{\nu}{2}})}$ and $\Gamma_2\subset \overline{(W_\eta\setminus W_{\frac{3\eta}{4}})}$. Therefore, for any $x\in \o_\nu\cap W_{2\eta}$
$$d(x,\Gamma_1 \cup \Gamma_2)\ge d\left(x,\overline{(\O(\o)\setminus \o_{\frac{\nu}{2}})}\cup\overline{(W_\eta\setminus W_{\frac{3\eta}{4}})} \right)>0.$$

\subsection*{Step 3:}
Let $\nu$ and $\eta$ be defined by the above steps.
Since $d>0$, choosing $\eps$ small enough, say $\eps\le \frac{d}{2\eta}$, it follows that 
$$\o_\nu\cap W_{2\eta}\subset \{ x|d\left(x,\partial(\O(\o)\cap W_\eta)\right)\ge \eps\eta\}.$$

Now, since $g\ge \eta$ in $\O(\o)\cap W_\eta$, from Lemma \ref{h-lem-loc-esti1} there exists $\eps^*$ so that for all $0<\eps\le\eps^*$ there exists a non empty set $\O_\eps$ and a constant $C_\eps(J,\eta,\O(\o)\cap W_\eta,b)$ such that
$$\left\{x|d(x,\partial(\O(\o)\cap W_\eta))\ge \eps\eta\right\} \subset \O_\eps\subset\left\{x|d(x,\partial(\O(\o)\cap W_\eta))\ge \frac{\eps\eta}{2}\right\}$$
and
\begin{equation}
\int_{\O_\eps}u(y)dy \ge C_\eps\int_{\O(\o)\cap W_\eta} u(y)dy\label{h-eq-sing-reg5}.
\end{equation} 

Observe that by choosing $\eps\le \min\{\eps^*,\frac{d}{2\eta}\}$, we also have $\o_\nu\cap W_{2\eta} \subset \O_\eps$.

We now fix $\eps\le \min\{\eps^*,\frac{d}{2\eta}\}$ and choose $\delta<\frac{\eps\eta}{8}$, and consider the set $$\O_\delta:=\bigcup_{x\in\O_\eps} B(x,2\delta).$$
 By construction, we have $\O_\eps\subset\subset\O_\delta\subset \O(\o)\cap W_\eta$ and $g\ge \eta$ in $\O_\delta$. Therefore, from Lemma \ref{h-lem-loc-esti2} there exists $\delta_1$ and $C'_0$ such that for any $x\in \O_\eps$, we have 
\begin{equation}
u(x)\ge C'_0 \int_{B(x,\delta_1)}u(s)ds.
\end{equation}
So we have 
\begin{equation}\label{h-eq-sing-reg6}
u(x)\ge \min\{C_0,C'_0\} \int_{B(x,\delta_1)}u(s)ds.
\end{equation}
\subsection*{Step 4:}
Take now $\delta^*\le \min\{\delta_0,\delta_1\}$, where $\delta_0$ is defined in \eqref{h-eq-sing-reg4}.
Then by construction $\O_\eps$ and $u$ satisfies the assumption of  Lemma \ref{h-lem-loc-esti3} and by (ii) of Lemma \ref{h-lem-loc-esti3} we end up with
%
\begin{equation}
 u(x)\ge \min\{C_0,C'_0\} \int_{B(x,\delta^*)}u(s)\,ds\ge \frac{C(\min\{C_0,C'_0\})}{N} \int_{\O_\eps}u(s)\,ds. \label{h-eq-sing-reg7}
\end{equation}
Collecting the inequalities \eqref{h-eq-sing-reg3}, \eqref{h-eq-sing-reg5} and \eqref{h-eq-sing-reg7}, it follows that 
 for all $x\in \O_\eps$ we have 
 
 \begin{align}
  u(x_0)&\le\frac{2\|J\|_{\infty}}{\inf_{\o}b(x) C_\eps}\int_{\O_\eps}u(s)\,ds\\
  & \le \frac{2\|J\|_{\infty} }{\inf_{\o}b(x) C_\eps}\frac{N}{C(\min\{C_0,C'_0\})} \int_{B(x,\delta^*)}u(s)\,ds\\
  &\le \frac{2\|J\|_{\infty} }{\inf_{\o}b(x) C_\eps}\frac{N}{C(\min\{C_0,C'_0\})\min\{C_0,C'_0\} }   u(x)\label{h-eq-sing-reg8}.
 \end{align}

Observing that  $\o_\nu\cap W_2\eta \subset \O_\eps$, from equation \eqref{h-eq-sing-reg8} it follows that 
\begin{equation}\label{h-eq-sing-reg9}
 u(x_0)\le  \bar C u(x) \qquad\text{ for all }\qquad x\in \o_\nu\cap W_2\eta.
\end{equation}

Now observe that from \eqref{h-eq-sing-reg2}, we also get that for all $x\in \o\setminus\o_\nu$, 
$$u(x)\ge C_0 \int_{B(\bar x,\delta_0)}u(s)\,ds, $$
where $B(\bar x, \delta_0)\subset \o_\nu\setminus \o_{4\nu}$.
Therefore, since $\bar x\in \o_\nu\cap W_{2\eta}\subset \O_\eps$ and $\delta^*\le \delta_0$, it follows that
\begin{equation}
u(x)\ge C_0 \int_{B(\bar x,\delta_0)}u(s)\,ds \ge \min\{C_0,C_0'\} \int_{B(\bar x,\frac{\delta^*}{2})}u(s)\,ds.\label{h-eq-sing-reg10}
\end{equation}

Using again Lemma \ref{h-lem-loc-esti3} yields
\begin{align}
 u(x)&\ge \frac{C( \min\{C_0,C_0'\})}{N} \int_{\O_\eps}u(s)\,ds .\label{h-eq-sing-reg11}
\end{align} 
As above, we can  combine the inequalities \eqref{h-eq-sing-reg11}, \eqref{h-eq-sing-reg5}  and \eqref{h-eq-sing-reg1} to obtain  
\begin{equation}\label{h-eq-sing-reg12}
 u(x_0)\le  \bar C u(x) \qquad\text{ for all }\qquad x\in (\o\setminus\o_\nu).
\end{equation}

Thus, from the inequalities  \eqref{h-eq-sing-reg9} and \eqref{h-eq-sing-reg12} we get 
\begin{equation}
 u(x_0)\le  \bar C u(y) \qquad\text{ for all }\qquad \, y \in \o\cap W_{2\eta}.
\end{equation}

Hence, we have 
\begin{equation}
 u(x)\le \bar C u(y) \qquad\text{ for all }\qquad x\in \o ,\, y \in \o\cap W_{2\eta}.
\end{equation}
\fdem

%

Let us now deal with the general estimate and prove the Theorem \ref{h-th2}.
  
\dem{Proof of Theorem \ref{h-th2}:}
Let us first observe that the estimate is straightforward if the set $\o\subset\subset \O$. Indeed,  in such  cases there exists always a regular compact set  $\tilde \o$ such that $\o \subset\subset\tilde \o\subset \O$ and the Theorem \ref{h-th3} applies.  Therefore the only case left to analyse is when the domain $\o$ touches the boundary of $\O$.
In such case and without any regularity assumption on the domain,  we cannot derive the estimates from a simple argument. To obtain the estimate we use a similar argument as for the proof of Theorem \ref{h-th3}.
The proof here follows essentially three steps.

\subsection*{Step 1:}
 First remark that as for the proof of Theorem \ref{h-th3}, we can restrict our attention to a compact set $\o\subset \O$ such that $\S\subset\subset \o$. 
Now let us define the sets $\O(\o):=\bigcup_{x\in \o}B(x,\beta)\cap \O$ and $W_\eta$ as in the above proof.
Following a similar argument, for a point $x_0\in \o$ where $u$ achieves its maximum and for small enough $\eta$, say $\eta\le \eta^*$, we have 
\begin{equation}
u(x_0)\le \frac{2\|J\|_{\infty}}{\inf_{\o}b(x)}\int_{\O(\o)\cap W_\eta }u(y) d\mu(y) \label{h-eq-sing-1}.
\end{equation}

From \eqref{h-eq}, using that $u$, $J$ and $g$ are non-negative, it follows that for all $x\in \O$ we have
\begin{equation}
u(x) \ge \frac{1}{b(x)}\int_{\O(\o)\cap W_\eta} J\left(\frac{x-y}{g(y)}\right) u(y) d\mu(y).\label{h-eq-sing-2}
\end{equation}

Since $g\ge \eta$ in $\O(\o)\cap W_\eta$, from Lemma \ref{h-lem-loc-esti1} there exists $\eps^*>0$ such that for all $0<\eps\le\eps^*$, there exists a non empty set $\O_\eps \subset\subset \O(\o)\cap W_\eta$ and a constant $C_\eps(J,\eta,\O(\o)\cap W_\eta,b)$ such that 
$$ \left\{x|d(x,\partial(\O(\o)\cap W_\eta))\ge \eps\eta\right\} \subset \O_\eps\subset\left\{x|d(x,\partial(\O(\o)\cap W_\eta))\ge \frac{\eps\eta}{2}\right\} $$
and 
\begin{equation}
\int_{\O_\eps}u(y)d\mu(y) \ge C_\eps\int_{\O(\o)\cap W_\eta} u(y)d\mu(y)\label{h-eq-sing-3}.
\end{equation} 
Recall that from the proof of Lemma \ref{h-lem-loc-esti1}, we also have $\eps\le \eps^*\le\frac{1}{4}$. Thus, we have 
$$ \left\{x|d(x,\partial(\O(\o)\cap W_\eta))\ge \eta\right\} \subset \left\{x|d(x,\partial(\O(\o)\cap W_\eta))\ge \eps\eta\right\}\subset \O_\eps.$$
\subsection*{Step 2:}
Choose now $\delta<\min\left\{\frac{\eps\eta}{8},\frac{\eps}{8}\right\}$, where $\eps$ and $\eta$ are defined by the previous Step, and consider the set $$\O_\delta:=\bigcup_{x\in\O_\eps} B(x,2\delta).$$
 By construction, $\O_\delta\subset \O(\o)\cap W_\eta$ and $g\ge \eta$ in $\O_\delta$. Using Lemma \ref{h-lem-loc-esti2}, for any $x\in \O_\eps$ we have 
\begin{equation}\label{h-eq-sing-4}
u(x)\ge C_0 \int_{B(x,\delta)}u(s)ds.
\end{equation}

\subsection*{Step 3}
%
As in the previous proof, by construction $\O_\eps$ and $u$ satisfies the assumption of  Lemma \ref{h-lem-loc-esti3} and we end up with
\begin{equation}
 u(x)\ge C_0 \int_{B(x,\delta)}u(s)\,ds\ge \frac{C(C_0)}{N} \int_{\O_\eps}u(s)\,ds. \label{h-eq-sing-5}
\end{equation}
Collecting \eqref{h-eq-sing-1}, \eqref{h-eq-sing-2} and \eqref{h-eq-sing-5}, it follows that 
 for all $x\in \O_\eps$ we have $$ u(x_0)\le  \frac{2\|J\|_{\infty} }{\inf_{\o}b(x) C_\eps}\frac{N}{C(C_0)C_0 } u(x).$$
Therefore, for all $y\in \o$ and $x \in \O_\eps$ we have 
 $$ u(y)\le \bar C u(x).$$

Hence, for all $y\in \o$ and $x \in \O_\eps\cap \o$ we have 
 $$ u(y)\le \bar C u(x).$$
\fdem

%
%


\section{Construction of non trivial positive solution of a particular nonlocal equation}
In this section, we deal with the construction of a positive solution of \eqref{h-eq-examp} and prove Theorem \ref{h-th4}. 

\dem{Proof of Theorem \ref{h-th4}:}
We treat two cases 
\subsection*{Case 1: $\O$ bounded}~\\
First, let us assume that $\O$ is bounded and $J$ is regular. Let us define the operator $T\in \l(C(\O))$ by 
$$Tu:=\frac{1}{a(x)}\int_{\O}J\left(\frac{x-y}{g(y)}\right)\frac{u(y)}{g^d(y)}\, dy,$$
where 
\begin{equation*}
a(x):=\left\{\begin{array}{l}
\int_{\O}J\left(\frac{y-x}{g(x)}\right)\frac{dy}{g^d(x)}\quad\text{for} \quad x\not \in \S\\
1\quad\text{otherwise}.
\end{array}
\right.
\end{equation*}
Since $\frac{1}{g^d(y)}\in L^{1}_{loc}$, $T$ is a compact operator. Moreover $T$ is positive since $g$,$a$ and $J$ are non-negative functions. Using now the Krein-Rutman Theorem, there exists an eigenvalue $\lambda$ and a continuous eigenfunction $\phi>0$ such that 
$$ \frac{1}{a(x)}\int_{\O}J\left(\frac{x-y}{g(y)}\right)\frac{\phi}{g(y)}\, dy=\lambda \phi.$$
Integrating the equation over $\O$, it follows that $\lambda=1$, and $\phi$ is our desired solution.
To obtain a solution in the case $J \in L^{\infty}\cap L^{1}(\O)$, we proceed by regularisation. 
Let $(\rho_n)_{n\in \N}$ a sequence of $C^\infty$ mollifier and consider the  problem with the regular kernel $J_n:=J\star\rho_n$. By the above argumentation,  for each $n$ there exists $\phi_n>0$ solution of the regularized problem.
Now using the global estimate  of  Corollary \eqref{h-cor2} we have 
$$sup_{\O} \phi_n\le C(J_n) \inf_{\o} \phi_n.$$
Let us normalize $\phi_n$ so that $\inf{\o} \phi_n \ge 1$.  From a carefully analysis of $C(J_n)$ on can show that this constant depends only on $C(J,\O)$.   Thus $(\phi_n)_{n\in \O}$ is a uniformly bounded sequence of continuous function. Using now Arzela-Ascoli Theoreme and a standard diagonal extraction argument, there exists a convergent subsequence $(\phi_n)$ which converges locally uniformly to a positive continuous bounded function $\phi$. Moreover, using Lebesgue monotone convergence  Theorem,  one can see that $\phi$ is our desired solution.

\subsection*{Case 2: $\O$ unbounded}~\\
Assume now that $\O$ is any open set and let $\O_n$ be an increasing sequence of bounded subsets such that $lim_{n\to \infty}\O_n=\O$.
Since $\S\subset\subset \O$, we can also assume that for all $n\in \N,\; \S\cap\O_n\subset\subset \O_n$.
Let $\phi_n$ denote the associated solution to $\O_n$ with the normalization $\phi(x_0)=1$ for some fixed $x_0 \in \O_n$ that we will choose later on. 
Since $n\in \N,\; \S\subset\subset \O_n$ and $(\O_n)_{n\in \N}$ is an increasing sequence of sets, for some $\eta_1$ small we have $$\bigcap_{n\in \N} (\O_n \cap W_{\eta_1})=\O_0\cap W_{\eta_1} \neq \emptyset. $$
Let us choose $x_0 \in \O_0 \cap W_{\eta_1}$.


Let us now fix $n\in \N$ and consider the sequence of functions $(\phi_{n+k})_{k\in\N}$.
By construction, $\phi_{n+k}$ satisfies the equation
$$\int_{\O_{n+k}}J\left(\frac{x-y}{g(y)}\right)\frac{\phi_{n+k}(y)}{g^d(y)}\, dy -a_{n+k}(x)\phi_{n+k}=0.$$
Since $\O_n$ is an increasing sequence of bounded sets, for any $k\in \N$ we have $\O_n\subset \O_{n+k}$.
Using Theorem \ref{h-th2} with $\O_n$ and $\phi_{n+k}$, it follows that for any $k\in \N$ there exists a constant $\eta_k^*$, such that for all $\eta\le \eta_k^*$ there exists $\o'_{n+k}$ and a constant $C_{n+k}(J,g,\eta,\|a_{n+k}\|_{\infty},\beta,\O_n)$ such that 
\begin{align}
&\{x\in \O_{n+k}|d(x,\partial(\O_n \cap W_\eta))>\eta\}\subset \o'_{n+k}\\
&\sup_{\O_n} \phi_{n+k}\le C_{n+k}(\eta) \phi_{n+k}(x) \quad\text{ for all } \quad x\in \o'_{n+k}\label{h-eq-sup-esti1}.
\end{align}
For each $k\in \N$, let us choose $\eta_k$ such that $$ \O_n\cap W_{\eta_1} \subset \{x\in \O_n|d(x,\partial W_{\eta_{k}})>\eta_k\}. $$
Using the monotonicity of the sequence $(\O_n)_{n\in\N}$, it follows that 
$$\O_0\cap W_{\eta_1} \subset \O_n \cap W_{\eta_1} \subset \{x\in \O_n|d(x,\partial W_{\eta_{k}})>\eta_k\}\subset \o'_{n+k}. $$
Therefore, from the above set inclusion and \eqref{h-eq-sup-esti1}, it follows that 
\begin{equation}\label{h-eq-sup-esti2}
\sup_{\O_n}\phi_n\le C_{n+k}(\eta_k) \phi_n(x_0)\le C_{n+k}(\eta_k).
\end{equation}

Now, observe that the sequence of positive functions $(a_{n+k}(x))_{k\in \N}$ is increasing in $\O_n$ and uniformly bounded. 
The monotonicity property follows easily from the monotonicity of the $\O_n$. Indeed, recall that for any $x \in \O_n\setminus \S$ we have
$$a_{n+k}(x)= \int_{\O_{n+k}}J\left(\frac{y-x}{g(x)}\right)\frac{dy}{g^d(x)}.$$
Therefore, using that $\O_n\subset \O_n+1$ and that $J,g$ are non negative functions, it follows that 
$$a_{n+k}(x)= \int_{\O_{n+k}}J\left(\frac{y-x}{g(x)}\right)\frac{dy}{g^d(x)}\le \int_{\O_{n+k+1}}J\left(\frac{y-x}{g(x)}\right)\frac{dy}{g^d(x)}=a_{n+k+1}(x).$$
On the other hand, for $x\in \S$, we have $a_n(x)=1$ for all $n\in \N$. Thus, $a_{n+k}\le a_{n+k+1}$ in $\O_n$. 

From the definition of $a_n$, we also get easily the uniform bound. 
For any $x\in \O_n\setminus \S$, using a change of variable we have
$$
a_{n}(x)= \int_{\O_{n}}J\left(\frac{y-x}{g(x)}\right)\frac{dy}{g^d(x)}\le\int_{\frac{\O_{n}-x}{g(x)}}J(z)\,dz\le 1.
$$

Using that $(a_n(x))_n$ is uniformly bounded independent of $n$ and increasing in $\O_n$, we can make the constant $C_{n+k}$ independent of $k$. 
Therefore, for all $k\in\N$, $\phi_{n+k}$ is uniformly bounded in $\O_n$. 
Now, since $\phi_{n+k}$ is uniformly continuous on $\O_n$, using Arzela-Ascoli's Theorem we can extract from $(\phi_{n+k})_{k\in\N}$ a sub-sequence which converges uniformly in $\O_n$. 
By a standard diagonal argument, we can extract from $(\phi_{n})_{n\in\N}$ a sub-sequence which converges to a function $\phi$ uniformly on every compact subset $\o$ of $\O$.
Using that $J$ has compact support and the Lebesgue dominated convergence theorem, passing to the limit in the equation yields 
$$\int_{\O}J\left(\frac{x-y}{g(y)}\right)\frac{\phi}{g(y)}\, dy-a(x) \phi=0.$$
 \fdem

\bigskip

\ackno{ The author thanks the Max Planck Institute for Mathematics in the Sciences for all the help they have provided during the realization of this work. The author would also warmly thanks the anonymous referee for his numerous pertinent comments which  have improved a lot the paper.  }

\bibliographystyle{plain}
\bibliography{harnack.bib}

\end{document}